\documentclass[10pt,reqno]{amsart}
\usepackage{graphicx, amssymb, color,pdfsync}
\usepackage[all,cmtip]{xy}
\numberwithin{equation}{section}
\usepackage{mathrsfs}
\usepackage{hyperref}

\hypersetup{
    colorlinks=true,       % false: boxed links; true: colored links
    linkcolor=blue,          % color of internal links
    citecolor=blue,        % color of links to bibliography
    filecolor=blue,      % color of file links
    urlcolor=blue           % color of external links
}
\usepackage{tikz}
\usetikzlibrary{matrix,arrows}
\DeclareMathOperator{\Aut}{Aut}
\usepackage[switch]{lineno}
\usepackage{yfonts}
\advance\hoffset by -0.9 truecm

\begin{document}
\newcommand{\s}{\vspace{0.2cm}}

\newtheorem{theo}{Theorem}
\newtheorem{prop}{Proposition}
\newtheorem{coro}{Corollary}
\newtheorem{lemm}{Lemma}
\newtheorem{claim}{Claim}
\newtheorem{example}{Example}
\theoremstyle{remark}
\newtheorem{rema}{\bf Remark}
\newtheorem{defi}{\bf Definition}

\title[A note on large automorphism groups of compact Riemann surfaces]{A note on large automorphism groups \\of compact Riemann surfaces}
\date{}

\author{Milagros Izquierdo}
\address{Matematiska institutionen, Link\"{o}pings Universitet}
\email{milagros.izquierdo@liu.se}

\author{Sebasti\'an Reyes-Carocca}
\address{Departamento de Matem\'atica y Estad\'istica, Universidad de La Frontera, Avenida Francisco Salazar 01145, Temuco, Chile.}
\email{sebastian.reyes@ufrontera.cl}

\thanks{Both authors were partially supported by Redes Etapa Inicial Grant 2017-170071. The second author was also partially supported by Fondecyt Grants 11180024, 1190991 and Anillo ACT1415 PIA-CONICYT Grant}
\keywords{Riemann surfaces, Group actions, Jacobian varieties}
\subjclass[2010]{14H30, 30F35, 14H37, 14H40}

\begin{abstract}  Belolipetsky and Jones classified those compact Riemann surfaces of genus $g$ admitting a large group of automorphisms of order $\lambda (g-1)$, for each $\lambda >6,$ under the assumption that $g-1$ is a prime number. In this article we study the remaining large cases; namely,  we classify Riemann surfaces admitting $5(g-1)$ and $6(g-1)$ automorphisms, with $g-1$ a prime number. As a consequence, we obtain the classification of Riemann surfaces admitting a group of automorphisms of order $3(g-1)$, with $g-1$ a prime number.  We also provide isogeny decompositions of their Jacobian varieties.  \end{abstract}
\maketitle

%%%%%%%%%%%%%%%%%%
%%%%%%%%%%%%%%%%%%
%%%%%%%%%%%%%%%%%%
%%%%%%%%%%%%%%%%%%
\section{Introduction and statement of the results}

The classification of groups of automorphisms of compact Riemann surfaces is a stimulating subject of study, and has attracted considerable interest since the nineteenth century. 

Let $S$ be a compact Riemann surface of genus $g \geqslant 2.$ It is well-known that the full automorphism group of $S$ is finite, and that its order is bounded by $84(g-1)$. 

A group of automorphisms $G$ of $S$ is said to be {\it large} if its order is strictly greater than $4(g-1);$ this bound arises naturally in the theory of Hurwitz spaces. In this case, it is known that $S$ is quasiplatonic (i.e.   cannot be deformed non-trivially in the moduli space together with its automorphisms) or belongs to a  complex one-dimensional family. See \cite{BCI,  CK, CI, K2}.

\s

Compact Riemann surfaces with large groups of automorphisms have been considered from different points of view. For instance, the cyclic case was considered by Wiman \cite{Wi}, Kulkarni  \cite{K2} and Singerman \cite{si3} (see also \cite{Hi}), and the abelian case was classified by Lomuto in \cite{AB}; see also \cite{PR}. Riemann surfaces with $8(g+1)$ automorphisms were considered by Accola \cite{Accola} and Maclachlan \cite{Mac}, and by Kulkarni in \cite{K1}. More recently, Riemann surfaces with $4g$ automorphisms were studied in \cite{BCI} (see also \cite{yojpaa}), and with $4(g+1)$ automorphisms in \cite{CI}. The maximal non-large case is considered in \cite{RR1}.

\s

Belolipetsky and Jones \cite{BJ} proved that under the assumption that $g-1$ is a prime number (sufficiently large to avoid sporadic cases), a compact Riemann surface of genus $g$ admitting a large group of automorphisms of order $\lambda(g-1)$, where $\lambda > 6,$ belongs to one of six infinite well-described sequences of Riemann surfaces. Similar results but stated in a combinatorial context of maps were obtained in \cite{pp}.

\s

In this article, we study and classify compact Riemann surfaces of genus $g \geqslant 8$ admitting a group of automorphisms of order $5(g-1)$ and $6(g-1),$ where $g-1$ a prime number; these cases were not considered in the article \cite{BJ} by Belolipetsky and Jones.  

We also determine an isogeny decomposition of the corresponding Jacobian varieties. 

The results of this paper are given in Theorem 1 and Theorem 2.

\s

\begin{theo} \label{t5q}
Let $g \geqslant 8$ such that $g-1$ is prime. There exists a compact Riemann surface $S$ of genus $g$ with a group of automorphisms $G$ of order $5(g-1)$ if and only if $g \equiv 2 \mbox{ mod }5.$ Moreover, in this case:

\begin{enumerate}
\item[(1)] the group $G$ is isomorphic to $$C_{g-1} \rtimes_5 C_5=\langle a,b : a^{g-1}=b^5=1, bab^{-1}=a^r\rangle,$$where $r$ is a primitive $5$-th root of unity in the field of $g-1$ elements, and $G$ acts with signature $(0; 5,5,5),$

\item[(2)] the action of $G$ extends to an action of a group $G'$ isomorphic to$$C_{g-1} \rtimes_{10}C_{10}=\langle a, c : a^{g-1} =c^{10}=1, cac^{-1}=a^{-r} \rangle,$$ with $r$ as before, and $G'$ acts with signature $(0; 2,5,10),$
\item[(3)] there are exactly four pairwise non-isomorphic such Riemann surfaces $S,$ 
\item[(4)] the full automorphism group of $S$ is $G',$ and
\item[(5)] the Jacobian variety $JS$ of each $S$ decomposes, up to isogeny, as the product $$JS \sim J(S/\langle a \rangle) \times (J(S/\langle c \rangle))^{10}.$$
  \end{enumerate}
\end{theo}

\begin{theo} \label{t6q}
Let $g \geqslant 8$ such that $g-1$ is prime.  There exists a compact Riemann surface of genus $g$ with a group of automorphisms of order $6(g-1)$ if and only if $g \equiv 2 \mbox{ mod }3.$ Moreover, in this case: 
\begin{enumerate}
\item[(1)] the Riemann surfaces form a closed one-dimensional equisymmetric family $\bar{\mathcal{F}}_{g}$ of Riemann surfaces $S$ with a group of automorphisms $G$ isomorphic to $$C_{g-1} \rtimes_6 C_6=\langle a,c : a^{g-1}=c^6=1, cac^{-1}=a^m\rangle,$$where $m$ is a primitive $6$-th root of unity in the field of $g-1$ elements, and $G$ acts with signature $(0; 2,2,3,3),$ 
\item[(2)] the Jacobian variety $JS$ of each $S$ in $\bar{\mathcal{F}}_{g}$ decomposes, up to isogeny, as the product $$JS \sim J(S/\langle a \rangle) \times (J(S/\langle c \rangle))^{6},$$
\item[(3)]  $\bar{\mathcal{F}}_g$ contains two Riemann surfaces $X_1$ and $X_2$ with a group of automorphisms $G'$ of order $12(g-1)$ isomorphic to $$(C_{g-1} \rtimes_6 C_6)\times C_2=\langle a,c, z : a^{g-1}=c^6=z^2=1, cac^{-1}=a^m, za=az, zc=cz \rangle $$with $m$ as before, acting with signature $(0;2,6,6),$ and  
\item[(4)] the Jacobian variety $JX_i$ of each $X_i$  can be decomposed, up to isogeny, as $$JX_i  \sim J(X_i/\langle a \rangle) \times (J(X_i/\langle cz \rangle))^6.$$
\end{enumerate}
Furthermore, if $\mathcal{F}_g$ denotes the interior of $\bar{\mathcal{F}}_g$ then:
\begin{enumerate}
\item[(5)] if $S \in \mathcal{F}_g$ then $G$ is the full automorphism group of $S,$ and
\item[(6)]  if $g > 14$ then the subset $\bar{\mathcal{F}}_g \setminus {\mathcal{F}}_g$ of $\bar{\mathcal{F}}_g$ is  $\{X_1, X_2\}, $ and the full automorphism group of $X_1$ and $X_2$ is $G'.$
\end{enumerate}
\end{theo}

\begin{rema} \mbox{}
\begin{enumerate}
\item The additional condition $g > 14$ in the last statement of Theorem \ref{t6q} is needed. In fact,  there are three non-isomorphic Riemann surfaces of genus 14 with full automorphism groups $H = \mbox{PSL}_2(13)$ of order $1092$ (these Riemann surfaces are Hurwitz curves and were considered by Macbeath in \cite{McB2}; see also \cite{C} and \cite{Streit}). Observe that the normaliser $G$ of a Sylow $13$-subgroup of $H$ is isomorphic to $C_{13} \rtimes_6 C_{6};$ thus, these Riemann surfaces belong to the closed family $\bar{\mathcal{F}}_{14}$ of Theorem \ref{t6q}. However, notice that $H$ does not contain any subgroup isomorphic to $(C_{13} \rtimes_{6} C_6)\times C_{2}$,  so these surfaces are not one of the surfaces $X_1, X_2$ in Theorem \ref{t6q}. 

\s

\item The family $\bar{\mathcal{F}}_g$ can be described from a hyperbolic geometry point of view by considering the deformations of a fundamental domain $\mathscr{P}$ of $\mathbb{H}$ with four non-equivalent vertices with total angles $\pi, \pi, 2\pi/3$ and $2\pi/3.$ The deformations of $\mathscr{P}$ correspond isometrically to the members of the family; see, for example, \cite{Farkas}, \cite{Nag} and \cite{singerman}. The two exceptional Riemann surfaces $X_1$ and $X_2$ arise geometrically when a deformation of $\mathscr{P}$ admits an extra involution.

\s

\item We anticipate the fact that whereas the Riemann surfaces in Theorem \ref{t5q} all appear in the article by Belolipetsky and Jones, the Riemann surfaces lying in the family $\bar{\mathcal{F}}_g$ of Theorem \ref{t6q} are  {\it new} (with only the exception of $X_1$ and $X_2$ which also appear in the article by Belolipetsky and Jones).
\end{enumerate}
\end{rema}

\s

As a consequence of the proof of the theorem above, we are able to easily derive a classification for the non-large case $\lambda=3.$

\begin{coro} \label{t3q}
Let $g \geqslant 8$ such that $g-1$ is prime. There exists a compact Riemann surface $S$ of genus $g$ with a group of automorphisms of order $3(g-1)$ if and only if $g \equiv 2 \mbox{ mod }3.$ Furthermore,  in this case 
 $S$ belongs to the family $\bar{\mathcal{F}}_g$ of Theorem \ref{t6q}. As a consequence, there is no compact Riemann surfaces of genus $g$ with full automorphism group of order 
 $3(g-1).$ 
\end{coro}

In Section \ref{preli} we will briefly review the background. The results will be proved in Sections \ref{prooft5q}, \ref{prooft6q} and \ref{prooft3q}. 

%; namely, Fuchsian groups, group actions on Riemann surfaces, the equisymmetric stratification of the moduli space, and the decomposition of Jacobian varieties with group action.

\section{Preliminaries} \label{preli}

\subsection{Fuchsian groups} Let $\mathbb{H}$ denote the upper-half plane, and let $\Gamma$ be a {\it cocompact Fuchsian group}; i.e.  a discrete group of automorphisms of $\mathbb{H}$ with compact orbit space $\mathbb{H}/\Gamma.$ The algebraic structure of $\Gamma$ is determined by its {\it signature}: \begin{equation} \label{sig} s(\Gamma)=(h; m_1, \ldots, m_l),\end{equation}where $h$ denotes the topological genus of the surface $\mathbb{H}/\Gamma,$ and $m_1, \ldots, m_l$ the branch indices in the (orbifold) universal covering $\mathbb{H} \to \mathbb{H}/\Gamma.$ If $l=0$, then $\Gamma$ is called a {\it surface Fuchsian group}. 

Each Fuchsian group $\Gamma$ has a {\it fundamental region}; namely, roughly speaking, a domain of $\mathbb{H}$ which contains exactly one point for each orbit of the action of $\Gamma$ on $\mathbb{H}$ (see, for example, \cite[Chapter IV]{Farkas}). 

Let $\Gamma$ be a Fuchsian group with signature \eqref{sig}. Then $\Gamma$ has a {\it canonical presentation}  with generators $a_1, \ldots, a_{h}$, $b_1, \ldots, b_{h},$ $ x_1, \ldots , x_l$ and relations
\begin{equation}\label{prese}x_1^{m_1}=\cdots =x_l^{m_l}=\Pi_{i=1}^{h}[a_i, b_i] \Pi_{i=1}^l x_i=1,\end{equation}where $[u,v]=uvu^{-1}v^{-1}.$ The hyperbolic area of each fundamental region of $\Gamma$ is given by $$\mu(\Gamma)=2 \pi [2h-2 + \Sigma_{j=1}^l(1-\tfrac{1}{m_j})]$$and the complex dimension of the Teichm\"{u}ller space associated to $\Gamma$ is $3g-3+l.$ See, for example, \cite{Farkas, Nag, singerman}.

Let $\Gamma'$ be a group of automorphisms of $\mathbb{H}.$ If a Fuchsian group $\Gamma$ is a finite index subgroup of $\Gamma'$ then $\Gamma'$ is also a Fuchsian group and their hyperbolic areas are related by the Riemann--Hurwitz formula $$\mu(\Gamma)= [\Gamma' : \Gamma] \cdot \mu(\Gamma').$$

\subsection{Riemann surfaces and group actions} \label{lapiz}

Let $S$ be a compact Riemann surface. We denote by $\mbox{Aut}(S)$ the full automorphism group of $S,$ and say that a group $G$ acts on $S$ if there is a group monomorphism $\psi: G\to \Aut(S).$ The space of orbits $S/G$ of the action of $G$ induced by $\psi(G)$ is naturally endowed with a Riemann surface structure such that the projection $S \to S/G$ is holomorphic.

By the uniformization theorem, a Riemann surface $S$ is conformally equivalent (isomorphic) to the quotient $\mathbb{H}/\Gamma$, where $\Gamma$ is a surface Fuchsian group. Lifting $G$ to the universal covering $\mathbb{H} \to \mathbb{H}/\Gamma,$ the group $G$ acts on  $S$ if and only if there is a Fuchsian group $\Gamma'$ containing $\Gamma$ and a group  epimorphism \begin{equation*}\label{epi}\theta: \Gamma' \to G \, \, \mbox{ such that }  \, \, \mbox{ker}(\theta)=\Gamma\end{equation*}(see \cite{Brou, Farkas, yoibero, singerman}). Such an epimorphism will be called {\em a surface epimorphism}. We say that the action of $G$ on $S$ is {\it given or represented} by the surface epimorphism $\theta.$ Note that the Riemann surface $S/G$ is isomorphic to $\mathbb{H}/\Gamma'$. We shall also say that $G$ acts on $S$ with signature $s(\Gamma').$

\s

Let us assume that $G$ is a subgroup of another group $G_1$. The action of $G$ on $S \cong \mathbb{H}/\Gamma$ is said to \textit{extend} to an action of $G_1$ if and only if there is a Fuchsian group $\Gamma''$ containing $\Gamma'$ together with a surface epimorphism $$\Theta: \Gamma'' \to G_1 \, \, \mbox{ in such a way that }  \, \, \Theta|_{\Gamma'}=\theta, \, \mbox{ker}(\Theta)=\mbox{ker}(\theta)=\Gamma,$$and $\Gamma'$ and $\Gamma''$ have associated Teichm\"{u}ller spaces of the same dimension. Singerman in \cite{singerman2} determined all those pairs of signatures $(s(\Gamma' ),  s(\Gamma''))$ for which it may be possible to extend actions. An action is called {\it maximal} if it cannot be extended in the aforementioned sense.

\subsection{Topologically equivalent actions} Let $S$ be a compact Riemann surfaces and let $\text{Hom}^+(S)$ denote the group of orientation preserving homeomorphisms of $S.$ Two actions $\psi_i: G \to \mbox{Aut}(S)$  are said to be {\it topologically equivalent} if there exist $\omega \in \Aut(G)$ and $h \in \text{Hom}^+(S)$ such that
\begin{equation}\label{equivalentactions}
\psi_2(g) = h \psi_1(\omega(g)) h^{-1} \hspace{0.5 cm} \mbox{for all} \,\, g\in G.
\end{equation}

Note that topologically equivalent actions have the same signature. Each orientation preserving homeomorphism $h$ satisfying (\ref{equivalentactions}) yields a group automorphism $h^*$ of $\Gamma'$ where $\mathbb{H}/\Gamma' \cong S/G$. We shall denote the subgroup of $\mbox{Aut}(\Gamma')$ consisting of the automorphisms $h^*$ by $\mathfrak{B}$. 
%The group $\mbox{Aut}(G) \times \mathscr{B}$ acts on the set of epimorphisms defining actions of $G$ on $S$ with signature $s(\Gamma')$ by $$((\omega, h^*), \theta) \mapsto \omega \circ \theta \circ (h^*)^{-1}.$$  

Two surface epimorphisms $\theta_1, \theta_2 : \Gamma' \to G$ define topologically equivalent actions if and only there are $\omega \in \mbox{Aut}(G)$ and $h^* \in \mathfrak{B}$ such that  $\theta_2 = \omega\circ\theta_1 \circ h^*$ (see \cite{Brou, Harvey, McB}). We remark that if the genus of $S/G$ is zero,  then the group $\mathfrak{B}$ is generated by the {\it braid transformations} $\Phi_{i,i+1} \in \mbox{Aut}(\Gamma')$ defined by: $$x_i \mapsto x_{i+1}, \hspace{0.3 cm}x_{i+1} \mapsto x_{i+1}^{-1}x_{i}x_{i+1} \hspace{0.3 cm} \mbox{ and }\hspace{0.3 cm} x_j \mapsto x_j \mbox{ when }j \neq i, i+1$$
for each $i \in \{1, \ldots, l-1\}.$ See, for example, \cite[p. 31]{trenzas} and also \cite{b, IJR}.

\subsection{Equisymmetric stratification} Let $\mathscr{M}_g$ denote the moduli space of compact Riemann surfaces of genus $g \geqslant 2.$ It is well-known that $\mathscr{M}_g$ is endowed with an orbifold structure and that its locus of orbifold-singular points, the so-called {\em branch locus} $\mathscr{B}_g,$  is formed by Riemann surfaces with non-trivial automorphisms for $g\geqslant 3$. For $g=2$ the branch locus $\mathscr{B}_2$ consists of the Riemann surfaces admitting other automorphisms than the hyperelliptic involution. See, for example, \cite{Nag}.

It was proved in \cite{b} (see also \cite{Harvey}) that the branch locus $\mathscr{B}_g$ admits an {\it equisymmetric stratification}, \begin{equation}\label{stratif} \mathscr{B}_g = \cup_{G, \theta} \bar{\mathscr{M}}_g^{G, \theta}\end{equation}where the non-empty {\it equisymmetric strata} are in bijective correspondence with the topological classes of  actions that are maximal (in the sense introduced in Subsection \ref{lapiz}). Concretely:

\begin{enumerate}

\item the {\it equisymmetric stratum} ${\mathscr{M}}_g^{G, \theta}$ consists of all those Riemann surfaces $S$ of genus $g$ with full automorphism group isomorphic to $G$ such that the action of $G$ on $S$ is topologically equivalent to $\theta$,

\item the closure $\bar{\mathscr{M}}_g^{G, \theta}$ of  ${\mathscr{M}}_g^{G, \theta}$ consists of all those Riemann surfaces $S$ of genus $g$ with a group of automorphisms isomorphic to $G$ such that the action of $G$ on $S$ is  topologically equivalent to $\theta$, 

\item $\bar{\mathscr{M}}_g^{G, \theta}$ is a closed irreducible algebraic subvariety of $\mathscr{M}_g$, and

\item  ${\mathscr{M}}_g^{G, \theta}$ is, if non-empty, a smooth, connected,
locally closed algebraic subvariety of $\mathscr{M}_{g}$ and Zariski dense in
$\bar{\mathscr{M}}_g^{G, \theta}$ 
\end{enumerate}

\s

In this work we shall use the following:

\s

{\bf Definition.} A closed family  $\bar{\mathcal{F}}$ of compact Riemann surfaces of genus $g$ whose members admit an action of a group $G$ will be called {\it equisymmetric} if its interior $\mathcal{F}$ consists of exactly one stratum.

\s

It is worth remarking that  the interior $\mathcal{F}$ of  $\bar{\mathcal{F}}$ consists of those Riemann surfaces whose full automorphism group is isomorphic to $G.$ Meanwhile, the subset $\bar{\mathcal{F}} \setminus \mathcal{F}$ corresponds to those Riemann surfaces that have strictly more automorphisms than $G.$

\subsection{Decomposition of Jacobian varieties} \label{jaco} Let $S$ be a compact Riemann surface of genus $g \geqslant 2.$ We denote by $JS$ the Jacobian variety of $S$, and recall that $JS$ is an irreducible principally polarized abelian variety of dimension $g.$ See, for example,  \cite{bl}.

The relevance of the Jacobian variety lies in the well-known Torelli's theorem, which asserts that two compact Riemann surfaces are isomorphic if and only if their Jacobian varieties are isomorphic as principally polarized abelian varieties.  

If a finite group $G$ acts on $S$ then this action induces an isogeny decomposition \begin{equation} \label{iso}JS \sim J(S/G) \times A_{2} \times \ldots \times A_{r}\end{equation} which is $G$-equivariant. The factors in \eqref{iso} are in bijective correspondence with the rational irreducible representations of $G$; the factor $A_1 \sim J(S/G)$ is associated to the trivial representation (see \cite{cr, l-r}). 

\s

The decomposition of Jacobian varieties with group actions has been extensively studied; the simplest case of such a decomposition was already noticed by Wirtinger in \cite{W} and used by Schottky and Jung in \cite{SJ}. For decompositions with respect to special groups, we refer to  \cite{d1, nos, IJR, PA, d3, RR2}.

\s
Let $G$ be a finite group. For each complex representation $\rho : G \to \mbox{GL}(V)$ of $G$ we shall denote its {\it degree} by $d_V$; i.e. the dimension of $V$ as a complex vector space. If $H$ is a subgroup of $G,$ then we shall denote the dimension of the vector subspace of $V$ fixed under the action $H$ by $d_V^H.$  By abuse of notation, we shall  write $V$ to refer to the representation $\rho.$ See \cite{Serre} for more details.

\s

Let us assume that $G$ acts on a Riemann surface $S$ with signature  \eqref{sig},  and that this action is determined by the surface epimorphism $\theta : \Gamma \to G.$ Let $H_1, \ldots, H_t$ be groups of automorphisms  of $S$ such that $G$ contains $H_i$ for each $i$. Following \cite{kanirubiyo}, the collection $\{H_1, \ldots, H_t\}$ is called {\it $G$-admissible} if $$d_{V}^{H_1}+ \cdots + d_{V}^{H_t}  \leqslant d_{V}$$
for every complex irreducible representation $V$ of $G$ in $\mathfrak{J},$ where the elements of $\mathfrak{J}$ are characterized (by using \cite[Theorem 5.12]{yoibero}) as follows:\begin{enumerate}
\item the trivial representation belongs to $\mathfrak{J}$ if and only if   the genus of $S/G$ is different from zero, and  
\item a non-trivial  representation $V$ belongs to $\mathfrak{J}$ if and only if $$d_{V}(\gamma -1)+\tfrac{1}{2}\Sigma_{i=1}^l (d_{V}-d_{V}^{\langle \theta(x_i) \rangle} ) \neq 0,$$
where the $x_i's$ are  canonical generators \eqref{prese} of $\Gamma.$
\end{enumerate}
The collection is called {\it admissible} if it is $G$-admissible for some group $G$.  The main result of \cite{kanirubiyo} ensures that if $\{H_1, \ldots, H_t\}$ is an admissible collection of  groups of automorphism of a Riemann surface $S$ then  $$JS \sim \Pi_{i=1}^t J(S/H_i) \times P$$
for some abelian subvariety $P$ of $JS.$ See also  \cite{KR}.

\s

{\bf Notation.} Let $n \geqslant 2$ be an integer and let $q$ be a prime.
Throughout this article we denote the cyclic group of order $n$  by $C_n$, the dihedral group of order $2n$ by $D_n$ and the field of $q$ elements by $\mathbb{F}_q$.

\section{Proof of Theorem \ref{t5q}}\label{prooft5q}Let $S$ be a compact Riemann surface of genus $g \geqslant 8,$ where $q=g-1$ is prime,  and assume that $S$ has a group of automorphisms $G$ of order $5q.$  If the signature of the action of $G$ is $(h; m_1, \ldots, m_l)$ then the Riemann--Hurwitz formula says that $$\tfrac{2}{5}=2h-2+l-\Sigma_{i=1}^l\tfrac{1}{m_i}.$$Note that if $h \geqslant 1$ then $l <1.$ Thus, $h=0$ and therefore $$\Sigma_{i=1}^l\tfrac{1}{m_i}=l-\tfrac{12}{5}.$$Now, the facts that each $m_i$ divides $5q$ and that $q \geqslant 7$ is prime imply that $l=3$ and $m_1=m_2=m_3=5.$ 

By the classical Sylow's theorems, if $q \not\equiv 1 \mbox{ mod }5$ then $G$ is isomorphic to $C_{5q},$ and if $q \equiv 1 \mbox{ mod }5$ then $G$ is isomorphic to either $C_{5q}$ or to $$C_{q}\rtimes_5 C_5= \langle a,b : a^q=b^5=1, bab^{-1}=a^r \rangle,$$where $r$ is a primitive $5$-th root of  unity in $\mathbb{F}_q$. Note that, since $C_{5q}$ cannot be generated by two elements of order five, if $q \not\equiv 1 \mbox{ mod }5,$ then 
there are no compact Riemann surfaces of genus $g$ with a group of automorphisms of order $5q.$

\s

From now on we assume that $q \equiv 1 \mbox{ mod }5$ and that $G \cong C_{q}\rtimes_5 C_5$.

\s

Let $\Gamma$ be a Fuchsian group of signature $(0; 5,5,5)$ with canonical presentation$$\Gamma=\langle x_1, x_2, x_3,: x_1^5=x_2^5=x_3^5=x_1x_2x_3=1\rangle$$and let $\theta: \Gamma  \to G$ be a surface epimorphism representing the action of $G$ on $S.$ We recall that $G$ has exactly four conjugacy classes of elements of order 5; namely $\{a^lb^j : 1 \leqslant l \leqslant q\} \mbox{ for } j=1,2,3,4.$

If the epimorphism $\theta$ is defined by $\theta(x_1)=a^{l_1}b^i, \, \theta(x_2)=a^{l_2}b^j  \mbox{ and } \theta(x_3)=a^{l_3}b^k$ where $l_1, l_2, l_3 \in \{1, \ldots, q\}$ and  $i,j,k \in \{1, \dots, 4\},$ then, after applying a suitable inner automorphism of $G$, we can assume $l_3 \equiv 0 \mbox{ mod } q$ and then $l_1\equiv -r^il_2 \mbox{ mod } q.$ As $l_2 \not\equiv 0 \mbox{ mod }q$ (otherwise $\theta$ is not surjective), we can consider the automorphism of $G$ given by $a \mapsto a^{t_2}$ and $b \mapsto b,$ where $l_2t_2 \equiv 1 \mbox{ mod } q,$ to see that $\theta$ is equivalent to the epimorphism $\theta_{i,j,k}$ defined by 
$$\theta_{i,j,k}(x_1)=a^{-r^i}b^i, \, \theta_{i,j,k}(x_2)=ab^j  \mbox{ and } \theta_{i,j,k}(x_3)=b^k.$$Now, as the braid automorphisms act by permuting conjugacy classes of elements of $\Gamma$ and as $i+j+k \equiv 0  \mbox{ mod } 5,$ there are at most four pairwise topologically non-equivalent actions of $G$ on $S$, represented by $$ \theta_1=\theta_{1,2,2}, \, \theta_2=\theta_{2,4,4}, \, \theta_3=\theta_{1,1,3} \mbox{ and } \theta_4=\theta_{3,3,4}.$$

Following \cite{singerman2}, the action given by each $\theta_{n}$ can be possibly extended to actions of signatures $(0; 3,3,5)$ and $(0; 2,5,10),$ and  these actions, in turn, can be possibly extended to a maximal action of signature $(0; 2,3,10).$  We claim that none of the actions of $G$ extends to an action of signature $(0;3,3,5).$ Indeed, if an action extends then there is a surface-kernel epimorphism $$\Delta(0; 3,3,5) \to H$$from a Fuchsian group of signature $(0; 3,3,5)$ onto a group $H$ of order $15q.$ However, by the Sylow's theorems, $H$ maps onto the cyclic group $C_{15}$ and therefore there is an epimorphism $$\Delta(0; 3,3,5) \to C_{15};$$
this is impossible. Note that this fact also ensures that none of the actions of $G$ extends to an action of signature $(0; 2,3,10).$ 

\s

Let us now consider a Fuchsian group $\Gamma_1$ of signature $(0;2,5,10)$ with canonical presentation 
$$\Gamma_1=\langle y_1, y_2, y_3 : y_1^2=y_2^5=y_3^{10}=y_1y_2y_3=1 \rangle,$$and a finite group $G' =C_{q} \rtimes_{10} C_{10}$ with presentation $$\langle a,b,s : a^{q}=b^5=s^2=1, bab^{-1}=a^r, sas=a^{-1}, [s,b]=1 \rangle,$$where $r$ is a primitive $5$-th root of unity in $\mathbb{F}_q$. As proved in \cite[Example (ii)]{BJ} (see also \cite[Theorem 3]{SW}), the surface epimorphisms $\Theta_n: \Gamma_1 \to G' \cong C_{q} \rtimes_{10} C_{10}$ given by $$\Theta_n(y_1)=as, \, \Theta_n(y_2)=ab^{2n} \mbox{ and } \Theta_n(y_3)=b^{-2n}s,$$ for $1 \leqslant n \leqslant 4,$ define four pairwise non-isomorphic Riemann surfaces $X_1, \ldots, X_4$ of genus $g$ with full automorphism group isomorphic to $C_{q} \rtimes_{10} C_{10}.$

 Note that the subgroup of $\Gamma_1$ generated by $$\tilde{x}_1=(y_1y_3)^{-1},   \, \tilde{x}_2 = y_2 \mbox{ and } \tilde{x}_3=y_3^2$$is isomorphic to $\Gamma,$ and that $\Theta_n(\tilde{x}_1)=a^{-r^{2n}}b^{2n},  \Theta_n(\tilde{x}_2)=ab^{2n}$ and  $\Theta_n(\tilde{x}_3)=b^{-4n}.$ It follows that $\Theta_n|_{\Gamma}=\theta_{n}\, \mbox{ for each } n \in \{1,2,3,4\}$ and therefore each action of $G \cong C_q \rtimes_5 C_5$ on $S$ with signature $(0; 5,5,5)$ extends to an action of $ G'\cong C_q \rtimes_{10}C_{10}$ with signature $(0; 2,5,10);$ thus $S$ is isomorphic to $X_i$ for some $i.$ In particular,  there does not exist a Riemann surface of genus $g$ with full automorphism group of order $5q.$

\s

Finally, we decompose the Jacobian variety $JS$ of each $S.$ If we set $c=bs$ and $m=-r$ then$$\mbox{Aut}(S)\cong \langle a,c : a^q=c^{10}=1, cac^{-1}=a^m\rangle=C_q \rtimes_{10} C_{10}.$$We shall use this presentation in the sequel. Set $ \omega_t:=\mbox{exp}(\tfrac{2 \pi i}{t}).$

\s

 The group $C_{q} \rtimes_{10} C_{10}$ has, up to equivalence, ten complex irreducible representations of  degree 1, given by $$U_i : a \mapsto 1, \, c \mapsto\omega_{10}^i$$for $0 \leqslant i \leqslant 9.$ Let $\alpha=\tfrac{q-1}{10} \in \mathbb{N}$ and choose integers $k_1, \ldots, k_{\alpha} \in \{1, \ldots, q-1\}$ in such a way that $$\sqcup_{j=1}^{\alpha} \{k_j, k_jm, k_jm^2, \ldots, k_jm^9\}=\{1, \ldots, q-1\},$$where $\sqcup$ stands for disjoint union. Then, the group $C_{q} \rtimes_{10} C_{10}$ has, up to equivalence,
$\alpha$ complex irreducible representations of degree 10, given by $$V_j: a \mapsto \mbox{diag}(\omega_q^{k_j}, \omega_q^{k_jm}, \omega_q^{k_jm^2}, \ldots, \omega_q^{k_jm^9}), \,\, c \mapsto \left( \begin{smallmatrix}
0 & 1 & 0 & \cdots & 0 \\
0 & 0 & 1 & \cdots & 0 \\
\, & \, & \, & \ddots & \, \\
0 & 0 & 0 & \cdots & 1 \\
1 & 0 & 0 & \cdots & 0 \\
\end{smallmatrix} \right)\mbox{ for } 1 \leqslant j \leqslant \alpha.
$$Consider $H=\langle a \rangle$ and $H_{t}=\langle a^{t}c \rangle$ for $1\leqslant t \leqslant 10,$ and notice that 
$$d_{U_i}^{H}+ \Sigma_{t=1}^{10} d_{U_i}^{H_{t}} = 1=d_{U_i} \,\, \mbox{ and } \,\, d_{V_j}^{H}+ \Sigma_{t=1}^{10} d_{V_j}^{H_{t}} = 10=d_{V_j}$$
for each $i \in \{1, \ldots, 9\}$ and for each $j \in \{1, \ldots, \alpha\}.$ Thereby, as explained in Subsection \ref{jaco}, the collection $\{H, H_1, \ldots, H_{10}\}$ is admissible and therefore, by \cite{kanirubiyo}, there is an abelian subvariety $P$ of $JS$ such that $$JS \sim J(S/H) \times \Pi_{t=1}^{10} J(S/H_t)  \times P\sim J(S/\langle a \rangle) \times (J(S/\langle c \rangle))^{10}  \times P,$$where the second isogeny follows after noticing that, for each $t,$ the groups  $H_t$ and $\langle c \rangle$ are conjugate.

Observe that the $q$-sheeted regular covering map $S \to S/\langle a \rangle$ is unbranched, and that the regular covering map $S \to S/\langle c \rangle$ ramifies over exactly three values, marked with 2, 5 and 10. Then, it follows from the Riemann--Hurwitz formula that the genera of $S/\langle a \rangle$ and $S/\langle c \rangle$ are 2 and $\alpha $ respectively; thus $P=0.$

\s
This completes the proof of Theorem \ref{t5q}.

\section{Proof of Theorem \ref{t6q}}\label{prooft6q}Let $S$ be a compact Riemann surface of genus $g \geqslant 8,$ where $q=g-1$ is prime,  and assume that $S$ has a group of automorphism $G$ of order $6q.$  Similarly as argued in the proof of Theorem \ref{t5q}, by the Riemann--Hurwitz formula the possible signatures of the action of $G$ on $S$ are $(0;2,2,3,3), (0;2,2,2,6)$ and $ (0;3,6,6)$ for each genus and, in addition, the signature $(0; 2,7,42)$ for $g=8.$  

First of all, the signature $(0; 2,7,42)$ for $g=8$ cannot be realized since there is no surface epimorphism from a Fuchsian group of signature $(0; 2,7,42)$ to a (necessarily cyclic) group of order 42. In addition, by the classical Sylow's theorems, $G$ contains exactly one normal subgroup isomorphic to $C_q$ and therefore $G$ is isomorphic to a semidirect product $C_q \rtimes H,$ where $H$ is a group of order 6.

\s
{\bf Claim 1.} $H=C_6.$
\s

Let us assume that $H=D_3$ and therefore$$ G \cong C_q \rtimes D_3= \langle a,b,s : a^q=b^3=s^2=1, (sb)^2=1, bab^{-1}=a^u, sas=a^{v} \rangle,$$where $u$ is either 1 or a primitive third root of unity in $\mathbb{F}_q$, and $v = \pm1.$ 

\begin{enumerate}
\item If $u=1$ and $v=1$ then $G$ is isomorphic to the direct product $C_q \times D_3.$ However, as among every collection of generators of $C_q \times D_3$ there must be an element of order a multiple of $q,$ we see that there are no compact Riemann surfaces of genus $g$ with a group of automorphisms isomorphic to $C_q \times D_3$ 
since the order of the generators of the Fuchsian groups corresponding to such actions are 2, 3 and 6.

\s

\item If $u=1$ and $v=-1$ then $G$ is isomorphic to  $D_{3q}$ and therefore $G$ has no elements of order 6.  Moreover,  the elements of order three are $(ab)^q$ and $(ab)^{2q}$, and the involutions are of the form $s(ab)^l$ for $1 \leqslant l \leqslant 3q.$ It can be checked that if the product of  two involutions and two elements of order three is 1, then these elements generate $D_6.$ All the above ensures that there are no Riemann surfaces of genus $g$ with a group of automorphisms isomorphic to $D_{3q}.$

\s

\item Finally, if $u$ is a primitive third root of unity in $\mathbb{F}_q$, then the equation $(sb)a(sb)^{-1}=a^{r^2v}$ shows that the action of the involution $sb$ on $C_q$ has order three for $v=1$ and order six for $v=-1$. This is not possible.
\end{enumerate}
This proves Claim 1.
\s
 
Thereby,  $G \cong C_q \rtimes C_6 = \langle a,b,s : a^q=b^3=s^2=1, [s, b]=1, bab^{-1}=a^u, sas=a^{v} \rangle$ where $u$ is either 1 or  a primitive third  root of unity in $\mathbb{F}_q$ and $v = \pm1.$ 

\s
{\bf Claim 2.} $u$ is a primitive third root of unity in $\mathbb{F}_q$.
\s

Assume $u=1.$ \begin{enumerate} 
\item If $v=1$, then $G \cong  C_{6q}$ which is not generated by elements of order two and three. Thus, there are no Riemann surfaces of genus $g$ with a group of automorphisms isomorphic to $C_{6q}$.

\s

\item If $v=-1$ then $G \cong C_q \rtimes_{2} C_6$ where $C_6$ acts on $C_q$ with order two.  The elements of order two are of the form $a^ls,$ the elements of order six of the form $a^lbs$ and $a^lb^2s$ for $1 \leqslant l \leqslant q,$ and the elements of order three are $b$ and $b^2$. It can be seen that: 
\begin{enumerate} 
\item $G$ cannot be generated by three elements, being two of them of order two and one of order three, in such a way that their product has order three, \item the product of three elements of order two must have order two, and
\item  $G$ cannot be generated by two elements of order six whose product has order three. 
\end{enumerate}
All the above ensures that there are no Riemann surfaces of genus $g$ with a group of  automorphisms isomorphic to $C_q \rtimes_{2} C_6.$ 
\end{enumerate}This proves Claim 2.

\s

Therefore, $G \cong C_q \rtimes C_6$ with a presentation 
$\langle a,b,s : a^q=b^3=s^2=1, [s, b]=1, bab^{-1}=a^r, sas=a^{v} \rangle$, where $v=\pm 1$ and $r$ is a primitive third root of unity in $\mathbb{F}_q$.  Consequently $g-1 = q \equiv 1 \mbox{ mod }3$.

\s
We have two cases for the finite group $G$:

\s
\noindent \textbf{Case 1.}
If $v=1$ then $G$ is isomorphic to $C_q \rtimes_{3} C_6$ where $C_6$ acts on $C_q$ with order three.  The elements of order 3 are of the form $a^lb$ and $a^lb^2,$ the elements of order 6 of the form $a^lbs$ and $a^lb^2s$ for $1 \leqslant l \leqslant q,$ and $s$ is the unique element of order two. 

\noindent \textbf{Case 2.}
If  $v=-1$ then $G$ is isomorphic to $C_q \rtimes_{6} C_6$ where $C_6$ acts on $C_q$ with order six. The elements of order two are of the form $a^ls,$ the elements of order three are of the form $a^lb$ and $a^lb^2,$ and the elements of order six of the form $a^lbs$ and $a^lb^2s$ for $1 \leqslant l \leqslant q.$

\s

We now study each possible signature separately.

\s

{\bf Signature $(0; 2,2,2,6).$} As in both groups $C_q \rtimes_{3} C_6$ and $C_q \rtimes_{6} C_6$ the product of three elements of order two has order two, we see that there is no group of order $6q$ acting on a Riemann surface of genus $g$ with signature $(0; 2,2,2,6).$ See also \cite{CI}

\s

{\bf Signature $(0;2,2,3,3).$}  We note that there are no compact Riemann surfaces of genus $g$ admitting an action of $C_q \rtimes_{3} C_6$ with signature $(0;2,2,3,3);$ this follows from the fact that $s$ (which is the unique involution) and an element of order three generate a group of order six. 

 By contrast, we show that there is a complex one-dimensional equisymmetric family $\bar{\mathcal{F}}_g$ of Riemann surfaces $S$ of genus $g$ with a group of automorphisms isomorphic to $C_q \rtimes_{6} C_6$ acting on $S$ with signature $(0;2,2,3,3).$ Indeed, let $\Gamma_3$ be a Fuchsian group of signature $(0;2,2,3,3)$ with canonical presentation$$\Gamma_3 = \langle x_1, x_2, x_3, x_4 = x_1^2=x_2^2=x_3^3=x_4^3=x_1x_2x_3x_4 =1\rangle.$$Then the surface epimorphism $\theta_{3,0} : \Gamma_3 \to C_q \rtimes_{6} C_6$ defined by$$\theta_{3,0}(x_1) = s, \, \, \theta_{3,0}(x_2) = as,  \, \, \theta_{3,0}(x_3) = ab^2 \, \mbox{ and }\, \theta_{3,0}(x_4) = b,$$provides the family $\bar{\mathcal{F}}_g$ of Riemann surfaces admitting an action of $C_q \rtimes_{6} C_6$. 

 To prove that $\bar{\mathcal{F}}_g$ is equisymmetric we notice that, up to a permutation of the generators of $\Gamma_3$, a surface epimorphism $\theta_3 : \Gamma_3 \to C_q \rtimes_{6} C_6$ is of the form$$\theta_3(x_1)= a^{l_1}s, \, \theta_3(x_2)=a^{l_2}s, \, \theta_3(x_3)= a^{l_3}b^2\, \mbox{ and }  \, \theta_3(x_4)= a^{l_4}b,$$for some $l_1, \ldots, l_4 \in \{1, \ldots, q\}.$ 
Moreover, after applying a suitable automorphism of $G$ of the form $a \mapsto a^u, b \mapsto a^vb$ we can suppose $l_1 \equiv 0 \mbox{ mod } q$ and $l_2\equiv 1 \mbox{ mod } q.$ Now, if we set $m=l_4$ then an epimorphism  $\theta_3$ is equivalent to one epimorphism $\theta_{3,m}$ given by $$\theta_{3,m}(x_1)=s,\, \theta_{3,m}(x_2)=as, \, \theta_{3,m}(x_3)=a^{1+(1+r)m}b^2, \, \theta_{3,m}(x_4)=a^mb, \,  1\leqslant m \leqslant  q$$As $\Phi_{3,4}^2 \cdot \theta_{3,m} =\Theta_{3,m+r-1},$  after iterating $\Phi_{3,4}^2$ a suitable number of times, we see that each epimorphism $\theta_{3,m}$ is equivalent to $\theta_{3,0}$, as desired. 

\s

 We claim that the full automorphism group of a Riemann surface in the interior $\mathcal{F}_g$ of $\bar{\mathcal{F}}_g$ is $G.$ Indeed, if all those Riemann surfaces lying in the interior of the family $\bar{\mathcal{F}}_g$ admit strictly more automorphisms, then by \cite{singerman2} the action could only extend to an action of a group of order $12q$ of signature $(0;2,2,2,3).$ The last situation is not possible by \cite[Theorem 2(a)]{BJ} for $q > 18$ and by \cite{C} for the remaining cases $q=7$ and $q=13$.

\s

{\bf Signature $(0; 3,6,6).$}  Let $\Gamma_1$ be a Fuchsian group of signature $(0; 3,6,6)$ and consider its canonical presentation$$\Gamma_1 = \langle x_1, x_2, x_3 = x_1^3=x_2^6=x_3^6=x_1x_2x_3 =1\rangle.$$Applying automorphisms of the finite group, we have that:
\begin{enumerate}
\item A surface epimorphism $\Gamma_1 \to C_q \rtimes_{3} C_6$ representing an action of $C_q \rtimes_{3} C_6$ on $S$ with signature $(0; 3,6,6)$ is equivalent to one defined by $$\theta_{1,i}(x_1)=b^i, \,\, \,\theta_{1,i}(x_2)=a^{-r^i}b^i s \, \mbox{ and } \, \theta_{1,i}(x_3)=a^ibs \hspace{0.5 cm} \mbox{ for } i=1 \mbox{ or } i=2.$$ 

\item A surface epimorphism $\Gamma_1 \to C_q \rtimes_{6} C_6$ representing an action of $C_q \rtimes_{6} C_6$ on $S$ with signature $(0; 3,6,6)$ is equivalent to the one defined by$$\theta_{2}(x_1)=ab, \,\,\, \theta_{2}(x_2)=b s \, \mbox{ and } \, \theta_{2}(x_3)=a^{r}bs.$$
\end{enumerate}

Using the results of \cite{singerman2}, we can ensure that the action of $G$ on $S$ can be extended possibly only to actions with signatures  $(0; 2,6,6)$ and $(0; 2,4,6).$

\s

Let $\Gamma_2$ be a Fuchsian group of signature $(0; 2,6,6)$ with canonical presentation$$\Gamma_2=\langle y_1, y_2, y_3 : y_1^2=y_2^6=y_3^6=y_1y_2y_3=1\rangle.$$Following \cite[Example (i)]{BJ}, there exist two non-isomorphic Riemann surfaces $X_1$ and $X_2$ of genus $g$ with a group of automorphisms of order $12q$ acting on them with signature $(0; 2,6,6).$ Furthermore,  $$\mbox{Aut}(X_i) \cong (C_q \rtimes_6 C_6) \times C_2 \, \mbox{ for } \, i=1,2,$$with corresponding non-equivalent surface epimorphisms $\Theta_i: \Gamma_2 \to (C_q \rtimes_6 C_6) \times C_2$ giving the actions of $\mbox{Aut}(X_i)$ on $X_i$ defined by:$$ \Theta_{1}(y_1)= as,  \,   \, \Theta_{1}(y_2)= bsz, \,   \Theta_{1}(y_3)=a^{-r^2}b^2z$$ $$\Theta_{2}(y_1)= as, \,   \, \Theta_{2}(y_2)= b^2sz, \,  \Theta_{2}(y_3)=a^{-r}bz,$$where $z$ generates the $C_2$ central factor. 

\s

{\bf Claim 3.} If $S$ is a compact Riemann surface with an action of a group of order $6q$ with signature $(0; 3,6,6)$ then $S$ is isomorphic to either $X_1$ or $X_2.$

\s
First of all, we have seen above that such an action is given by the surface epimorphisms $\theta_{1,i}$ and $\theta_2.$ We see now that these actions extend. Setting $x_1'=y_2^2, x_2'=y_3$ and $x_3'=(y_2^2y_3)^{-1}$, the subgroup of $\Gamma_2$ generated by $x_1', x_2', x_3'$ is isomorphic to $\Gamma_1.$ Moreover  $$\Theta_1(x_1')= b^2,  \,   \, \Theta_1(x_2')= a^{-r^2}b^2z, \,  \, \Theta_1(x_3')=ab^2z,$$   $$\Theta_2(x_1')= b,  \,   \, \Theta_2(x_2')= a^{-r}bz, \,  \, \Theta_2(x_3')=abz.$$

 Note that $\langle a,b,z \rangle \cong C_q \rtimes_3 C_6$ and that the restrictions$$\Theta_1|_{\langle x_1', x_2', x_3'\rangle}, \Theta_2|_{\langle x_1', x_2', x_3'\rangle} : \Gamma_1 \cong \langle x_1', x_2', x_3'\rangle \to  C_q \rtimes_3 C_6$$are precisely $\theta_{1,2}$ and $\theta_{1,1}$ respectively. It follows that the action $\theta_{1,1}$  and  $\theta_{1,2}$ of $C_q \rtimes_{3} C_6$ on  compact Riemann surfaces $S$ of genus $g$ with signature $(0; 3,6,6)$ extend to the action of $(C_q \rtimes_{6} C_6) \times C_2$ with signature $(0; 2,6,6)$ represented by $\Theta_2$ and $\Theta_1$ respectively;  thus, $S$ is isomorphic to $X_2$ in the first case, and $S$ is isomorphic to $X_1$ in the second case. 
 
\s

Now, setting $x_1''=y_3^2, x_2''=y_2$ and $x_3''=(y_3^2y_2)^{-1}$, the subgroup of $\Gamma_2$ generated by $x_1'', x_2'', x_3''$ is isomorphic to $\Gamma_1.$ Moreover$$\Theta_1(x_1'')= ab,  \,  \Theta_1(x_2'')= b(sz), \,  \Theta_1(x_3'')=a^{-r}b(sz),$$ $$ \Theta_2(x_1'')= ab^2,     \, \Theta_2(x_2'')= b^2(sz), \, \Theta_2(x_3'')=a^{-r^2}b^2(sz).$$Note that $\langle a,b,sz \rangle \cong C_q \rtimes_6 C_6$ and that the restrictions$$\Theta_1|_{\langle x_1'', x_2'', x_3''\rangle}, \Theta_2|_{\langle x_1'', x_2'', x_3''\rangle} : \Gamma_1 \cong \langle x_1'', x_2'', x_3''\rangle \to  C_q \rtimes_6 C_6$$are equivalent to $\theta_2.$ It follows that the action $\theta_2$ of $C_q \rtimes_{6} C_6$ on a Riemann surface $S$ with signature $(0; 3,6,6)$ extends to both actions of $(C_q \rtimes_{6} C_6) \times C_2$ with signature $(0; 2,6,6)$ represented by $\Theta_1$ or by $\Theta_2;$ thus, $S$ is isomorphic to $X_1$ in the first case, and is isomorphic to $X_2$ in the second case. 
 
 This proves Claim 3.
 
 \s
 
Note that $ \hat{x}_1= (y_1y_2^2y_3^2)^{-1}, \hat{x}_2=y_1, \hat{x}_3= y_2^2$ and $\hat{x}_4= y_3^2$ generate a subgroup $\hat{\Gamma}$ of $\Gamma_2$ isomorphic to a Fuchsian group of signature $(0; 2,2,3,3).$ Furthermore, the restrictions $\Theta_1|_{\hat{\Gamma}}$ and $\Theta_2|_{\hat{\Gamma}}$ are epimorphisms equivalent to $\theta_{3,0}.$ This yields that $X_1$ and $X_2$ lie in $\bar{\mathcal{F}}_g \setminus \mathcal{F}_g$ as desired.

\s 
 
Finally,  we apply \cite[Theorem 2(a)]{BJ}  to conclude that for $g > 14$:

 \begin{enumerate} 
 \item the Riemann surfaces $X_1$ and $X_2$ are the unique compact Riemann surfaces with a group of automorphisms of order $12q,$ 
 \item there are no compact Riemann surfaces of genus $g$ with $24q$ automorphisms (in particular, the action of $G$ on $S$ of signature $(0; 3,6,6)$ cannot be extended to an action of signature $(0; 2,4,6)$),  and therefore $\{ X_1, X_2\} = \bar{\mathcal{F}}_g \setminus \mathcal{F}_g.$
 \end{enumerate}
 
 \s

We now decompose the associated Jacobian varieties; to do that we proceed  analogously as done in the proof of Theorem \ref{t5q}. Let $S \in \bar{\mathcal{F}}_g$ and set $\omega_t:=\mbox{exp}(\tfrac{2 \pi i}{t}).$  

Note that the group$$\mbox{Aut}(S) \cong C_{q} \rtimes_{6} C_{6}=\langle a,c : a^q=c^{6}=1, cac^{-1}=a^n\rangle$$where $n$ is a primitive $6$-th root of unity in $\mathbb{F}_q$,  has, up to equivalence, six complex irreducible representations of  degree 1, given by $$U_i : a \mapsto 1, \, c \mapsto\omega_{6}^i$$for $0 \leqslant i \leqslant 5.$ In addition, $C_{q} \rtimes_{6} C_{6}$ has  
$\beta = \tfrac{q-1}{6} \in \mathbb{N}$  complex irreducible representations of degree 6, namely$$V_j: a \mapsto \mbox{diag}(\omega_q^{k_j}, \omega_q^{k_jn}, \ldots, \omega_q^{k_jn^5}), \,\, c \mapsto \left( \begin{smallmatrix}
0 & 1 & 0 & 0 & 0 & 0 \\
0 & 0 & 1 & 0 & 0  & 0 \\
0 & 0 & 0 & 1 & 0  & 0 \\
0 & 0 & 0 & 0 & 1  & 0 \\
0 & 0 & 0 & 0 & 0  & 1 \\
1 & 0 & 0 & 0 & 0  & 0 \\
\end{smallmatrix} \right)\mbox{ for } 1 \leqslant j \leqslant \beta,
$$where $k_1, \ldots, k_{\beta} \in \{1, \ldots, q-1\}$ are integers chosen to satisfy that $$\sqcup_{j=1}^{\beta} \{k_j, k_jn, k_jn^2, \ldots, k_jn^5\}=\{1, \ldots, q-1\},$$where $\sqcup$ denotes the disjoint union. 

\s
Consider the subgroups $H=\langle a \rangle$ and $H_{t}=\langle a^{t}c \rangle$ for $t \in \{1, \ldots, 6\}.$ Note that $$d_{U_i}^{H}+ \Sigma_{t=1}^{6} d_{U_i}^{H_{t}} = 1=d_{U_i} \, \mbox{ and } \, d_{V_j}^{H}+ \Sigma_{t=1}^{6} d_{V_j}^{H_{t}} = 6=d_{V_j}$$for each $i \in \{1, \ldots, 6\}$ and for each $j \in \{1, \ldots, \beta\}.$  Thereby,  the collection $\{H, H_1, \ldots, H_{6}\}$ is admissible and therefore, by \cite{kanirubiyo}, there is an abelian subvariety $Q$ of $JS$ such that $$JS \sim J(S/H) \times \Pi_{t=1}^{6} J(S/H_t)  \times Q\sim J(S/\langle a \rangle) \times (J(S/\langle c \rangle))^{6}  \times Q,$$where the second isogeny follows from the fact that, for each $t,$ the groups $H_t$ and $\langle c \rangle$ are conjugate. 

\s
The $q$-sheeted regular covering map $S \to S/\langle a \rangle$ is unbranched, and the regular covering map $S \to S/\langle c \rangle$ ramifies over exactly four values, two marked with 2 and two marked with 3. Thus, the Riemann--Hurwitz formula implies that the genera of $S/\langle a \rangle$ and $S/\langle c \rangle$ are $2$ and $\beta$ respectively; thus $Q=0.$

\s

Let $S$ be one of the two non-isomorphic Riemann surfaces with $12q$ automorphisms. Each complex irreducible representation of $\mbox{Aut}(S) \cong (C_{q} \rtimes_{6} C_{6}) \times C_2=\langle a,c \rangle \times \langle z\rangle$ coincides with the tensor product of a complex irreducible representation of  $C_{q} \rtimes_{6} C_{6}$ and one of $C_2$ (see, for example \cite[p. 27]{Serre}) Thus, keeping the same notations as above, we see that the complex irreducible representations of  $(C_{q} \rtimes_{6} C_{6}) \times C_2$ are $$U_i^{\pm} : a \mapsto 1, \, c \mapsto\omega_{6}^i, \,z \mapsto \pm 1$$ for each $0 \leqslant i \leqslant 5,$ and $$V_j^{\pm}: a \mapsto \mbox{diag}(\omega_q^{k_j}, \omega_q^{k_jn}, \ldots, \omega_q^{k_jn^5}), \,\, c \mapsto \left( \begin{smallmatrix}
0 & 1 & 0 & 0 & 0 & 0 \\
0 & 0 & 1 & 0 & 0  & 0 \\
0 & 0 & 0 & 1 & 0  & 0 \\
0 & 0 & 0 & 0 & 1  & 0 \\
0 & 0 & 0 & 0 & 0  & 1 \\
1 & 0 & 0 & 0 & 0  & 0 \\
\end{smallmatrix} \right), \, z \mapsto \pm I_6$$for $1 \leqslant j \leqslant \beta,$ where $I_6$ denotes the $6 \times 6$ identity matrix.

\s

If we write $N=\langle a \rangle$ and $N_{t}=\langle a^{t}cz \rangle$ for $t \in \{1, \ldots, 6\},$ then it can be checked that the collection $\{N, N_1, \ldots, N_{6}\}$ is admissible. In addition,  as $N_t$ and $\langle cz \rangle$ are conjugate, we apply the result of  \cite{kanirubiyo} to ensure the existence of an abelian subvariety $R$ of $JS$ such that $$JS \sim J(S/\langle a \rangle) \times (J(S/\langle cz \rangle))^{6}  \times R.$$

The $q$-sheeted regular covering map $S \to S/\langle a \rangle$ is unbranched and  the regular covering map $S \to S/\langle cz \rangle$ ramifies over exactly three values, two marked with 2 and one with marked 3.  The Riemann--Hurwitz formula implies that the genera of $S/\langle a \rangle$ and $S/\langle cz \rangle$ are $2$ and $\beta$ respectively; thus $R=0.$

This finishes the proof of Theorem \ref{t6q}.

\section{Proof of Corollary \ref{t3q}} \label{prooft3q}
Let $S$ be a compact Riemann surface of genus $g \geqslant 8,$ where $q=g-1$ is prime,  and assume that $S$ has a group of automorphisms $G$ of order $3q.$ Similarly as argued in the proof of Theorem \ref{t5q},  by the Riemann--Hurwitz formula the possible signatures for the action of $G$ on $S$ are $(1;3)$ and $(0; 3,3,3,3)$ for each $g$ and, in addition, the signature $(0; 7,7,21)$ for $g=8.$ The latter exceptional case for $g=8$ can be disregarded because there are no surface epimorphisms from a Fuchsian group of signature $(0; 7,7,21)$ to a (necessarily cyclic) group of order 21.

\s
By the classical Sylow's theorems if $q \not\equiv 1 \mbox{ mod }3$ then $G$ is isomorphic to $C_{3q},$ and if $q \equiv 1 \mbox{ mod }3$ then $G$ is isomorphic to either $C_{3q}$ or to $$C_{q}\rtimes_3 C_3 = \langle a,b : a^q=b^3=1, bab^{-1}=a^r \rangle,$$ where $r$ is a primitive third root of unity in $\mathbb{F}_q$.    

As $C_{3q}$ is abelian, and as the commutator subgroup of $C_{q}\rtimes_3 C_3$ does not have elements of order three, we see that there are no compact Riemann surfaces of genus $g$ with a group of automorphisms of order $3q$ acting with signature $(1;3).$ Furthermore, as  $C_{3q}$ cannot be generated by elements of order three, we obtain that if $q \not\equiv 1 \mbox{ mod }3$ then there are no compact Riemann surfaces of genus $g$ with a group of automorphisms of order $3q$ acting with signature $(0;3,3,3,3).$ 

\s

Thus, from now on we assume that $g-1 = q \equiv 1 \mbox{ mod }3$ and that $G \cong C_{q}\rtimes_3 C_3$.

\s

Let $\Gamma'$ be a Fuchsian group of signature $(0; 3,3,3,3)$ with canonical presentation$$\Gamma'=\langle x_1, x_2, x_3,x_4 : x_1^3=x_2^3=x_3^3=x_4^3=x_1x_2x_3x_4=1\rangle$$and let $\theta: \Gamma'  \to G$ be a surface epimorphism representing the action of $G$ on $S.$ We recall that $G$ has exactly two conjugacy classes of elements of order 3: $\mathscr{C}_1=\{a^lb : 1 \leqslant l \leqslant q\} \mbox{ and } \mathscr{C}_2=\{a^lb^2 : 1 \leqslant l \leqslant q\}.$

\s
Note that among the elements $\theta(x_1), \ldots, \theta(x_4)$ of $G$ exactly two of them must belong to $\mathscr{C}_1$; otherwise their product is different from 1. Up to a permutation we can suppose that $\theta(x_i)=a^{l_i}b^2$ for $i=1,2$ and $\theta(x_i)=a^{l_i}b$ for $i=3,4,$ for suitable $l_1, \dots, l_4.$  Note that if $l_1 \equiv l_2 \mbox{ mod } q$ and $l_3 \equiv l_4 \mbox{ mod } q,$ then $\theta$ is not surjective; thus, after  considering an automorphism of $G$ of the form $a \mapsto a^i, b \mapsto a^jb,$ we have two cases to consider:

\begin{enumerate}
\item[(1)] Case A: $l_1 \equiv 0 \mbox{ mod } q$ and $l_2 \equiv 1 \mbox{ mod } q.$ 
\item[(2)] Case B: $l_3 \equiv 0 \mbox{ mod } q$ and $l_4 \equiv 1 \mbox{ mod } q.$ 
\end{enumerate}

\s

By \cite{singerman2}, the action of $G$ on $S$ can possibly be extended to an action of signature $(0; 2,2,3,3).$ We shall prove that --in both cases-- the action does extend and that the surfaces $S$ belong to the family $\bar{\mathcal{F}}_g$ of Theorem \ref{t6q}. To accomplish this task, let us consider the Fuchsian  group $\Gamma_3$ of signature $(0; 2,2,3,3)$ with canonical presentation$$\Gamma_3=\langle y_1, y_2, y_3, y_4 : y_1^2=y_2^2=y_3^{3}=y_1^3=y_1y_2y_3y_4=1 \rangle$$and let $G' \cong C_q \rtimes_{6} C_6$ be the finite group with presentation$$ \langle a,b,s : a^{q}=b^3=s^2=1, bab^{-1}=a^r, sas=a^{-1},  [s,b]=1\rangle,$$as in Section \ref{prooft6q}. We recall that, following the proof of Theorem \ref{t6q}, each surface epimorphism $\Theta : \Gamma_3 \to  C_q \rtimes_{6} C_6$ representing an action of $C_q \rtimes_{6} C_6$ is equivalent to the epimorphism $\Theta_m$ $(= \theta_{3,m}$ in the notation of the proof of Theorem \ref{t6q}) given by$$\Theta_m(y_1)=s,\, \Theta_m(y_2)=as, \, \Theta_m(y_3)=a^{1+(1+r)m}b^2, \, \Theta_m(y_4)=a^mb$$for a suitable $1 \leqslant m \leqslant q.$ 

\s

We now study each case separately.

\s

{\bf Case A.} In this case $\theta$ is equivalent to the epimorphism $\theta_n$ defined by$$\theta_n(x_1)=b^2, \,\, \theta_n(x_2)=ab^2, \,\, \theta_n(x_3)=a^{-r(n+1)}b, \, \,\theta_n(x_4)=a^nb, \, \mbox{ with } \, 1 \leqslant n \leqslant q.$$We notice that  
the subgroup of $\Gamma_3$ generated by $$\tilde{x}_1=y_3, \,\, \tilde{x}_2 = y_4^2, \,\, \tilde{x}_3=y_4 \,\mbox{ and }\, \tilde{x}_4=y_3^2$$is isomorphic to $\Gamma'$ and $$\Theta_m(\tilde{x}_1)=a^{1+(1+r)m}b^2, \, \, \Theta_m(\tilde{x}_2)=a^{(1+r)m}b^2, \, \, \Theta_m(\tilde{x}_3)=a^{m}b, \,\, \Theta_m(\tilde{x}_4)=a^{m-r}b.$$ Now, after considering the automorphism of $G$ given by $a \mapsto a^{-1}, b \mapsto ab^j \, \mbox{ where } \, j=\tfrac{1+m(1+r)}{1-r},$ it follows that $\Theta_m|_{\Gamma'}=\theta_{\frac{2(mr+r+1)}{1-r}} \, \mbox{ for each } \, m \in \{1, \ldots, q\}.$ Thereby, 
 each action of type A of $C_{q}\rtimes_3 C_3$ with signature $(0; 3,3,3,3)$ extends to an action of $C_q \rtimes_{6} C_6$ with signature $(0; 2,2,3,3).$ 

\s

 {\bf Case B.} In this case $\theta$ is equivalent to the epimorphism $\theta_l$ defined by$$\theta_l(x_1)=a^lb^2, \,\, \theta_l(x_2)=a^{-1-lr}b^2, \,\, \theta_l(x_3)=b, \, \,\theta_l(x_4)=ab, \, \mbox{ with } \, 1 \leqslant l \leqslant q.$$

Note that the subgroup of $\Gamma_3$ generated by $$\tilde{x}_1:=(y_2y_1)^{l}((y_2y_1)^{-m}y_4)^2, \,\,\, \tilde{x}_2:=(y_2y_1)^{-1-lr}((y_2y_1)^{-m}y_4)^5$$ $$\tilde{x}_3:=(y_2y_1)^{-m}y_4, \,\mbox{ and }\, \tilde{x}_4:=(x_1x_2x_3)^{-1}$$is isomorphic to $\Gamma'$ and $$\Theta_m(\tilde{x}_1)=a^lb^2, \,\, \Theta_m(\tilde{x}_2)=a^{-1-lr}b^2, \,\, \Theta_m(\tilde{x}_3)=b, \,\, \Theta_m(\tilde{x}_4)=ab.$$

Thereby,  
 each action of type B of $C_{q}\rtimes_3 C_3$ with signature $(0; 3,3,3,3)$ extends to an action of $C_q \rtimes_{6} C_6$ with signature $(0; 2,2,3,3).$

\s

As a consequence, there do not exist compact Riemann surfaces with full automorphism group of order $3q,$ and the proof of Corollary \ref{t3q} is complete.

\s

{\bf Remark 2.} The Riemann surfaces of genus $g=3$ admitting the action of a group of order six or twelve are given and classified in \cite{Brou}. There is no Riemann surface of genus three admitting an automorphism of order five. The Riemann surfaces of genus $g=4$ with a group of automorphisms of order fifteen or eighteen are given and classified in \cite{bci} and  \cite{CI2}.  Among them there is the equisymmetric family of cyclic trigonal surfaces with two trigonal morphisms (see \cite{CIY2, CIY, gabino}). Finally, the Riemann surfaces of genus $g=6$ with twenty-five or thirty automorphisms are given in \cite{magaard}.

\s

{\bf Acknowledgments.} The authors are grateful to Gareth Jones and to the referee for their helpful comments and suggestions. This article was mainly written when the second author visited Link\"{o}ping University; he wishes to express his gratitude for the hospitality and kindness of the Department of Mathematics during his research stay there.


\begin{thebibliography}{999}
\bibitem{Accola}
{\sc R. Accola,} {\em On the number of automorphisms of a closed Riemann surface}, Trans. Am. Math. Soc., {\bf 131} (1968), 398--408.

\bibitem{bci} 
{\sc G. Bartolini, A. F. Costa, and M. Izquierdo,} 
{\em On the orbifold structure of the moduli space of Riemann surfaces of genera four and five.}
Rev. R. Acad. Cienc. Exactas Fis. Nat. Ser. A Mat. RACSAM {\bf 108} (2014), no. 2, 769--793. 

\bibitem{BJ} 
{\sc M. V. Belolipetsky and G. A. Jones,} 
{\em Automorphism groups of Riemann surfaces of genus $p + 1$, where $p$ is prime.}
Glasg. Math. J. {\bf 47} (2005), no. 2, 379--393.

\bibitem{bl}
{\sc Ch. Birkenhake and H. Lange},
 { \em Complex Abelian Varieties,} $2^{nd}$ edition,
Grundl. Math. Wiss. {\bf 302}, Springer, 2004.


%\bibitem{Bo}
%{\sc{O. V. Bogopol'skii,}}
%{\em{Classifying the actions of finite groups on orientable surfaces of genus 4,}} [translation of Proceedings of the Institute of Mathematics, 30 (Russian), 48--69, Izdat. Ross. Akad. Nauk, Sibirsk. Otdel., Inst. Mat., Novosibirsk, 1996]. 
%Siberian Advances in Mathematics. 
%Siberian Adv. Math. {\bf 7} (1997), no. 4, 9--38.

%\bibitem{Br}
%{\sc{Th. Breuer},} 
%{\em{Characters and automorphism groups of compact Riemann surfaces}}, London Math. Soc. Lect. Notes {\bf 280}, Cambridge Univ. Press 2000.

\bibitem{Brou}
{\sc{S. A. Broughton,}}  {\em{Finite groups actions on surfaces of low genus,}} J. Pure Appl. Algebra {\bf 69} (1991), no. 3, 233--270.

\bibitem{b}
{\sc S. A. Broughton},
 { \em The equisymmetric stratification of the moduli space and the Krull dimension of mapping class groups,} Topology Appl. {\bf 37} (1990), no. 2, 101--113.
 
\bibitem{BCI}
{\sc{E. Bujalance, A. F. Costa and M. Izquierdo,}}  {\em{On Riemann surfaces of genus g with 4g automorphisms,}} Topology Appl. {\bf 218} (2017) 1--18.

\bibitem{cr}
{\sc A. Carocca and R. E. Rodr\'iguez,}
{\em Jacobians with group actions and rational idempotents.}
J. Algebra \textbf{306} (2006), no. 2, 322--343.

\bibitem{d1}
{\sc A. Carocca, S. Recillas and R. E. Rodr\'iguez},
 { \em Dihedral groups acting on Jacobians,} Contemp. Math. {\bf 311} (2011), 41--77.

\bibitem{C}
{\sc M. Conder,} www.math.auckland.ac.nz/$\sim$conder/OrientableRegularMaps101.txt 

\bibitem{CK}
{\sc M. Conder and R. Kulkarni,} 
{\em Infinite families of automorphism groups of Riemann surfaces}, in Discrete Groups and Geometry, Birmingham, 1991, in: Lond. Math. Soc. Lect. Note Ser., vol. 173, Cambridge Univ. Press. Cambridge, 1992, pp. 47--56.


\bibitem{pp}
{\sc M. Conder, J. \v{S}ir\'a\v{n} and T.
Tucker}, 
{\em The genera, reflexibility and simplicity of regular maps}.  J. Eur. Math. Soc. (JEMS) {\bf 12} (2010), no. 2, 343--364.


\bibitem{CI}
{\sc A. F. Costa and M. Izquierdo,}
{\em One-dimensional families of Riemann surfaces of genus $g$ with $4g + 4$ automorphisms,} 
Rev. R. Acad. Cienc. Exactas Fis. Nat. Ser. A Mat. RACSAM {\bf 112} (2018), no. 3, 623--631.

\bibitem{CI2}
{\sc A. F. Costa and M. Izquierdo,}
{\em Equisymmetric stratification of the singular locus of the moduli space of Riemann surfaces of genus 4.} Lecture Notes London Mathematical Society {\bf 368} Cambridge University Press 2010, 120-138



\bibitem {CIY2}
{\sc A. F. Costa, M. Izquierdo, and D. Ying,} 
{\em On cyclic $p$-gonal Riemann surfaces with several $p$-gonal morphisms}. Geom. Dedicata {\bf 147}
(2010), 139-147

\bibitem {CIY}
{\sc A. F. Costa, M. Izquierdo, and D. Ying,} 
{\em On trigonal Riemann
surfaces with non-unique morphisms}. Manuscripta Mathematica, {\bf 118}
(2005), 443-453

\bibitem{Farkas}
{\sc H. Farkas and I. Kra}, {\em Riemann surfaces}, Grad. Texts in Maths. {\bf 71}, Springer-Verlag (1980).

\bibitem{gabino}
{\sc G. Gonz\'alez-Diez,} {\em On prime Galois coverings of the Riemann sphere.} Ann. Mat. Pura Appl.  {\bf 168} (1995), 1--15. 

\bibitem{Harvey}
{\sc J. Harvey,}
{\em On branch loci in Teichm\"{u}ller space},
Trans. Amer. Math. Soc. {\bf 153} (1971), 387--399.

%\bibitem{Harvey1}
%{\sc J. Harvey,}
%{\em Cyclic groups of automorphisms of a compact Riemann surface}. 
%Quarterly J. Math. {\bf 17}, (1966), 86--97.

\bibitem{nos}
{\sc R. A. Hidalgo, L. Jim\'enez, S. Quispe and S. Reyes-Carocca,} {\em Quasiplatonic curves with symmetry group $\mathbb{Z}_2^2 \rtimes \mathbb{Z}_m$ are definable over $\mathbb{Q}$,} Bull. London Math. Soc. {\bf 49} (2017) 165--183.

%\bibitem{hur}
%{\sc A. Hurwitz,} {\em \"{Uber algebraische Gebilde mit eindeutigen Transformationen in sich,}} Math. Ann. {\bf 41} (1893), 403--442.

\bibitem {Hi}
{\sc S. Hirose,} {\em On periodic maps over surfaces with large periods}.
Tohoku Math. J. (2) \textbf{62} (2010), no. 1, 45--53.


\bibitem{IJR}
{\sc M. Izquierdo, L. Jim\'enez, A. Rojas,}
{\em Decomposition of Jacobian varieties of curves with dihedral actions via equisymmetric stratification},   Rev. Mat. Iberoam. {\bf 35}, No. 4 (2019), 1259--1279.

\bibitem{KR}
{\sc E. Kani and M. Rosen,} {\em Idempotent relations and factors of Jacobians},
Math. Ann. {\bf 284} (1989), 307--327.

%\bibitem{Ku}
%{\sc A. Kuribayashi and I. Kuribayashi,} {\em Automorphism groups of compact Riemann surfaces of genus three and four},
%J. Pure Appl. Algebra {\bf 65} (1990), 277--292.

\bibitem{trenzas}
{\sc C. Kassel and V. Turaev,} {\em Braid Groups},
Graduate Texts in Mathematics. {\bf 247} Springer-Verlag, New York.

\bibitem{K1}
{\sc R. S. Kulkarni,}
{ \em A note on Wiman and Accola-Maclachlan surfaces}.
Ann. Acad. Sci. Fenn., Ser. A 1 Math. {\bf 16} (1) (1991)
83--94.

\bibitem{K2}
{\sc R. S. Kulkarni,}
{ \em Riemann surfaces admitting large automorphism groups}, in: Extremal Riemann Surfaces, San Francisco, CA, 1995, in: Contemp. Math., vol. 201, Amer. Math. Soc., Providence, RI, 1997, pp. 63--79.

%\bibitem{Kis}
%{\sc R. S. Kulkarni,}
%{ \em Infinite families of surface symmetries}, 
%Israel J. Math. {\bf 76} (1991), no. 3, 337--343. 


\bibitem{l-r}
{\sc H. Lange and S. Recillas,}
{ \em Abelian varieties with group actions}.
J. Reine Angew. Mathematik, \textbf{575} (2004) 135--155.

\bibitem{AB}
{\sc C. Lomuto,}
{ \em Riemann surfaces with a large abelian group of automorphisms},
Collect. Math. {\bf 57} (2006), no. 3, 309--318. 

\bibitem{Mac}
{\sc C. Maclachlan,} 
{\em A bound for the number of automorphisms of a compact Riemann surface},
J. London Math. Soc. {\bf 44} (1969), 265--272.



\bibitem{McB2}
{\sc A. MacBeath,}
{\em Generators of the linear fractional groups}. 1969
Number Theory (Proc. Sympos. Pure Math., Vol. XII, Houston, Tex., 1967) 14--32 Amer. Math. Soc., Providence, R.I.



\bibitem{McB}
{\sc A. MacBeath,}
{\em The classification of non-euclidean crystallographic groups,} Canad. J. Math. {\bf 19} (1966), 1192--1205.

\bibitem{magaard}
{\sc  K. Magaard, T. Shaska, S. Shpectorov and H. V\"olklein,} 
{\em  The locus of curves with prescribed automorphism group}.
Comm. Arithm. Fund. Groups (Kyoto),
Surikaisekikenkyusho Kokyuroku No. \textbf{1267} (2002), 112--141.

\bibitem{Nag}
{\sc{S. Nag,}} {\em{The complex analytic theory of Teichm\"{u}ller spaces,}} Canadian Math. Soc. Series of Monographs and Advanced Texts, Wiley-Intersciences (1988).

\bibitem{PA}
{\sc J. Paulhus and A. M. Rojas,}
{ \em Completely decomposable Jacobian varieties in new genera},  Experimental Mathematics {\bf 26} (2017), no. 4, 430--445.

\bibitem{PR}
{\sc R. Pignatelli and C. Raso,}
{ \em Riemann surfaces with a quasi large abelian group of automorphisms},  Matematiche (Catania) {\bf 66} (2011), no. 2, 77--90. 


\bibitem{d3}
{\sc S. Recillas and R. E. Rodr\'iguez,} {\em Jacobians and representations of $S_3$}, Aportaciones Mat. Investig. {\bf 13}, Soc. Mat. Mexicana, M\'exico, 1998.

\bibitem{yojpaa}
{\sc S. Reyes-Carocca,} {\em On the one-dimensional family of Riemann surfaces of genus $q$ with $4q$ automorphisms}, J. of Pure and Appl. Algebra 223, no. {\bf 5} (2019), 2123--2144.

\bibitem{RR1}
{\sc S. Reyes-Carocca,} {\em On Riemann surfaces of genus $g$ with $4g-4$ automorphisms}, To appear in: Israel Journal of Mathematics (2019), 	arXiv:1812.00705.

\bibitem{kanirubiyo}
{\sc S. Reyes-Carocca and R. E. Rodr\'iguez,} {\em A generalisation of Kani-Rosen decomposition theorem for Jacobian varieties}, Ann. Sc. Norm. Super. Pisa Cl. Sci. (5) {\bf 19} (2019), no. 2, 705--722.

\bibitem{RR2}
{\sc S. Reyes-Carocca and R. E. Rodr\'iguez,} {\em On Jacobians with group action and coverings} To appear in: Mathematische Zeitschrift (2019), DOI :10.1007/s00209-019-02263-3, arXiv:1711.07552v1 

\bibitem{yoibero}
{\sc A. M.   Rojas}, {\em Group actions on Jacobian varieties}, Rev. Mat. Iber. {\bf 23} (2007), 397--420.

\bibitem{Serre}
{\sc J. P. Serre}, {\em Linear Representations of Finite Groups}, Graduate Texts in Mathematics. {\bf 42} Springer-Verlag, New York.

\bibitem{SJ}
{\sc F. Schottky and H. Jung}, {\em Neue S\"atze \"uber Symmetralfunctionen und due Abel'schen Functionen der Riemann'schen Theorie,} S.B. Akad. Wiss. (Berlin) Phys. Math. Kl. {\bf 1} (1909), 282--297.

\bibitem{singerman2}
{\sc D. Singerman}, {\em Finitely maximal Fuchsian groups}, J. London Math. Soc. (2)  {\bf 6}, (1972), 29--38.

\bibitem{singerman}
{\sc D. Singerman}, {\em Subgroups of Fuchsian groups and finite permutation groups}, Bull. London Math. Soc.  {\bf 2}, (1970), 319--323.

\bibitem{si3}
{\sc D. Singerman}, {\em Symmetries of Riemann surfaces with large automorphism group}, Math. Ann.  {\bf 210}, (1974), 17--32.

\bibitem{SW}
{\sc M. Streit and J. Wolfart}, {\em Characters and Galois invariants of regular dessins,}  Rev. Mat. Complut. {\bf 13} (2000), no. 1, 49--81.

\bibitem{Streit}
{\sc M. Streit}, {\em Field of definition and Galois orbits for the Macbeath-Hurwitz curves,}  Arch. Math. (Basel) {\bf 74} (2000),  342--349.


\bibitem {Wi}
{\sc  A. Wiman}, 
{\em  \"{U}ber die hyperelliptischen Curven und diejenigen von Geschlechte p - Jwelche eindeutige Tiansformationen in sich zulassen.} \textit{Bihang till K. Svenska Vet.-Akad. Handlingar, Stockholm} \textbf{21}
(1895-6) 1-28.


\bibitem{W}
{\sc W. Wirtinger}, {\em Untersuchungen \"uber Theta Funktionen,} Teubner, Berlin, 1895.

\end{thebibliography}
\end{document}